\documentclass{article}
\usepackage{amsmath}
\usepackage{amssymb}
\usepackage{amsfonts}
\usepackage{amsthm}

\def\reals{\mathbb{R}}

\def\uball{\mathbb{B}}
\def\ereals{\overline{\mathbb{R}}}

\def\comp{\raise 1pt \hbox{$\scriptstyle\circ$}}

\def\argmin{\mathop{\rm argmin}\limits}

\def\minimize{\mathop{\rm minimize}\limits}
\def\maximize{\mathop{\rm maximize}\limits}

\def\st{\mathop{\rm subject\ to}}

\def\dom{\mathop{\rm dom}\nolimits}

\def\lev{\mathop{\rm lev}\nolimits}

\def\upto{{\raise 1pt \hbox{$\scriptstyle \,\nearrow\,$}}}
\def\downto{{\raise 1pt \hbox{$\scriptstyle \,\searrow\,$}}}

\def\co{\mathop{\rm co}}

\def\FF{(\F_t)_{t=0}^T}

\def\A{{\cal A}}
\def\B{{\cal B}}
\def\C{{\cal C}}
\def\D{{\cal D}}

\def\F{{\cal F}}
\def\G{{\cal G}}

\def\L{{\cal L}}

\def\N{{\cal N}}

\def\U{{\cal U}}

\def\Y{{\cal Y}}

\newtheorem{theorem}{Theorem}
\newtheorem{lemma}[theorem]{Lemma}
\newtheorem{corollary}[theorem]{Corollary}

\newtheorem{example}{Example}
\newtheorem{remark}{Remark}
\theoremstyle{definition}

\title{Stochastic programs without duality gaps}
\author{Teemu Pennanen \and Ari-Pekka Perkki\"o}

\begin{document}

\maketitle

\begin{abstract}
This paper studies dynamic stochastic optimization problems parametrized by a random variable. Such problems arise in many applications in operations research and mathematical finance. We give sufficient conditions for the existence of solutions and the absence of a duality gap. Our proof uses extended dynamic programming equations, whose validity is established under new relaxed conditions that generalize certain no-arbitrage conditions from mathematical finance.
\end{abstract}

\section{Introduction}

Let $(\Omega,\F,P)$ be a probability space with a filtration $\FF$ (an increasing sequence of sub-sigma-algebras of $\F$) and consider the dynamic stochastic optimization problem
\begin{equation}\label{p}\tag{P}
\minimize\quad Ef(x(\omega),u(\omega),\omega)\quad\text{over $x\in\N$},
\end{equation}
where, for given integers $n_t$ and $m$
\[
\N = \{(x_t)_{t=0}^T\,|\,x_t\in L^0(\Omega,\F_t,P;\reals^{n_t})\},
\]
$u\in L^0(\Omega,\F,P;\reals^m)$ and $f$ is an extended real-valued convex {\em normal integrand} on $\reals^n\times\reals^m\times\Omega$, where $n=n_0+ \ldots +n_T$. Recall that $L^0(\Omega,\F_t,P;\reals^{n_t})$ denotes the space of equivalence classes of $\F_t$-measurable $\reals^{n_t}$-valued functions that coincide $P$-almost surely. That $f$ is a normal integrand, means that the set-valued mapping $\omega\mapsto\{(x,u,\alpha)\,|\,f(x,u,\omega)\le\alpha\}$ is closed-valued and $\F$-measurable; see e.g.\ \cite[Chapter~14]{rw98}. This implies that $f$ is $\F\otimes\B(\reals^n\times\reals^m)$-measurable (see \cite[Corollary~14.34]{rw98}), so that $\omega\mapsto f(x(\omega),u(\omega),\omega)$ is $\F$-measurable for every $x\in\N$. Throughout this paper, the expectation is defined for any measurable function by setting it equal to $+\infty$ unless the positive part is integrable. We will also assume that $\F$ as well as $\F_t$ for $t=0,\ldots,T$ are complete with respect to $P$\footnote{This allows us to use certain results on conditional expectations of integrands which are not necessarily normal in the general case. This is based on \cite[Corollary~14.34]{rw98}, which says that, when $\F$ is $P$-complete, then a function $f:\Omega\times\reals^n\times\reals^m\to\ereals$ is a normal integrand if and only if it is $\F\otimes\B(\reals^n\times\reals^m)$-measurable and $(x,u)\mapsto f(\omega,x,u)$ is lower semicontinuous for every $\omega$.}.

The measurable function $u$ may be interpreted as a parameter or a perturbation in a given stochastic optimization problem. It was shown in \cite{pen11c} that \eqref{p} covers many important problems in operations research and mathematical finance and how the conjugate duality framework of Rockafellar~\cite{roc74} allows for a unified treatment of many well-known duality frameworks. In that context, the lower-semicontinuity of the value function
\[
\varphi(u)=\inf_{x\in\N}Ef(x(\omega),u(\omega),\omega)
\]
over an appropriate space of measurable functions $u$ is equivalent to the absence of a duality gap; see \cite[Section~2]{pen11c} for a precise statement. In certain applications, most notably in mathematical finance, the objective in \eqref{p} lacks the inf-compactness properties required by the classical ``direct method'' of calculus of variations for establishing lower semicontinuity (and the existence of solutions). It was shown in \cite[Section~5]{pen11c} how certain measure theoretic techniques from mathematical finance can be combined with classical techniques of convex analysis to obtain the lower semicontinuity of $\varphi$. It is essential for this that the strategies $x\in\N$ are allowed to be general measurable functions not restricted to be e.g.\ integrable. The lower semicontinuity result given in \cite{pen11c}, however, applies to normal integrands $f$ that take only the values $0$ and $+\infty$. While that already covers some fundamental results in mathematical finance, as illustrated in \cite[Section~6]{pen11c}, it is far from satisfactory from the general point of view.

The main purpose of this paper is to establish the lower semicontinuity of $\varphi$ for more general normal integrands. This will be done in Section~\ref{sec:lsc}. Our proof extends that of \cite[Theorem~8]{pen11c}, which employs a recursive argument reminiscent of dynamic programming. We clarify this connection in Section~\ref{sec:dp} by generalizing the dynamic programming equations proposed by Rockafellar and Wets~\cite{rw76} for stochastic convex optimization. The dynamic programming equations were substantially generalized already by Evstigneev~\cite{evs76} who removed many of the assumptions made in \cite{rw76}, including convexity. We will show that in the convex case, the inf-compactness assumption made in both \cite{rw76} and \cite{evs76} can be replaced by weaker ``recession condition'' which subsumes, in particular, various no-arbitrage conditions used in mathematical finance. An early application of recession analysis to utility maximization in financial markets can be found in Bertsekas~\cite{ber74}. Section~\ref{sec:app} of this paper gives an application to an optimal consumption problem in illiquid markets.

\section{Dynamic programming}\label{sec:dp}

The purpose of this section is to extend the {\em dynamic programming} recursion of \cite[Section~3]{rw76} which generalizes the classical Bellman equation for convex stochastic optimization. We will use the notion of a conditional expectation of a normal integrand much as in \cite{evs76} where certain assumptions (convexity, nonanticipativity of the domain of $f$ and the boundedness of the strategies) of \cite{rw76} were relaxed. We show that, in the convex case, the inf-compactness assumption used in both \cite{rw76} and \cite{evs76} can be replaced by a milder condition on the directions of recession much like in the classical closedness results of finite-dimensional convex analysis; see \cite[Section~8]{roc70a}. In certain financial applications, the new condition turns out to be equivalent to the classical no-arbitrage condition.

Let $X$ be a nonnegative $\F$-measurable function and let $\G\subseteq\F$ be another sigma-algebra. Then, there is a $\G$-measurable nonnegative function $E^\G X$, unique up to sets of $P$-measure zero, such that
\begin{equation}\label{eq:ce}
E[\chi_AX]=E[\chi_A(E^\G X)]\quad\forall A\in\G,
\end{equation}
where $\chi_A$ denotes the characteristic function of $A$; see e.g.\ Shiryaev~\cite[II.7]{shi96}.
The function $E^\G X$ is called the $\G$-{\em conditional expectation} of $X$. For a general $\F$-measurable extended real-valued function $X$, we set
\[
E^\G X:=E^\G X_+-E^\G X_-,
\]
where again, the convention $\infty-\infty=\infty$ is used. It is easily checked that with the extended definition of the integral, \eqref{eq:ce} is then valid for any measurable function $X$. Our definition of conditional expectation extends \cite[Definition~II.7.1]{shi96}, which assumes that $\min\{E^\G X_+,E^\G X_-\}<\infty$ almost surely. Our choice of setting $\infty-\infty=\infty$ is not arbitrary but specifically directed towards minimization problems.

The $\G$-{\em conditional expectation} of a normal integrand $h$ is a $\G$-measurable normal integrand $E^\G h$ such that
\[
(E^\G h)(x(\omega),\omega)= E^\G [h(x(\cdot),\cdot)](\omega)\quad P\text{-a.s.}
\]
for all $x\in L^0(\Omega,\G,P;\reals^n)$. 
There are various conditions that guarantee the existence and uniqueness of a conditional expectation of a normal integrand; see e.g.\ Bismut~\cite{bis73}, Dynkin and Evstigneev~\cite{de76}, Castaing and Valadier~\cite[Section~VIII.9]{cv77}, Thibault~\cite{thi81}, Truffert~\cite{tru91} or Choirat, Hess and Seri~\cite{chs3}. The following suffices for the purposes of this paper.

\begin{lemma}\label{lem:ce}
Let $\G\subseteq\F$ be a sigma-algebra and assume that $h$ is an $\F$-normal integrand with an integrable lower bound i.e.\ an integrable function $m$ such that $h(x,\omega)\ge m(\omega)$ for every $x$ and $\omega$. Then $h$ has a well-defined conditional expectation $E^\G h$ which has the integrable lower bound $E^G m$.
\end{lemma}

\begin{proof}
The integrable lower bound implies, for example, the quasi-integrability condition of Thibault~\cite{thi81} as well as the condition of Choirat, Hess and Seri~\cite{chs3}, both of which give the existence and uniqueness of the conditional expectation. It follows from the monotonicity of the conditional expectation that if $f\ge m$ for an integrable function $m$, then $E^\G f\ge E^\G m$.
\end{proof}

We will study problem \eqref{p} for a fixed $u\in L^0(\omega,\F,P;\reals^m)$ so we will omit it from the notation and define 
\[
h(x,\omega)=f(x,u(\omega),\omega).
\]
By \cite[14.45(c)]{rw98}, $h$ is a normal integrand. The convexity of $f$ implies that of $h$. We will use the notation $E_t=E^{\F_t}$ and $x^t=(x_0,\ldots,x_t)$ and define extended real-valued functions $h_t,\tilde h_t:\reals^{n_1+\dots+n_t}\times\Omega\rightarrow\ereals$ recursively for $t=T,\ldots,0$ by
\begin{equation}\label{dp}
\begin{split}
\tilde h_T&=h,\\
h_t &= E_t\tilde h_t,\\
\tilde h_{t-1}(x^{t-1},\omega)&=\inf_{x_t\in\reals^{n_t}}h_t(x^{t-1},x_t,\omega).
\end{split}
\end{equation}
This is essentially the dynamic programming recursion introduced in \cite{rw76}. Our formulation with conditional expectations of normal integrands is closer to \cite{evs76}, where certain assumptions of \cite{rw76} were relaxed. In the above formulation, one does not separate the decision variables $x_t$ into  ``state'' and ``control'' like in the classical dynamic programming models; see e.g.\ \cite{ber76} and \cite{bs78}. A formulation closer to the classical dynamic programming equations will be given in Corollary~\ref{cor:bel} below. A recent application of dynamic programming to mathematical finance can be found in R{\'a}sonyi and Stettner~\cite[Section~5]{rs5}.

In order to ensure that $h_t$ and $\tilde h_t$ are well-defined it suffices to require that the function $h$ has an integrable lower bound and that $h(\cdot,\omega)$ is inf-compact (i.e.\ $\{x\in\reals^n\,|\,h(x,\omega)\le \alpha\}$ is compact for every $\alpha\in\reals$) for every $\omega\in\Omega$; see \cite[Theorem~5]{evs76}. In the convex case, the compactness assumption can be replaced by a weaker condition stated in terms of the {\em recession function} of $h$. If $\dom h(\cdot,\omega)$ is nonempty, then the recession function has the expression
\[
h^\infty(x,\omega)=\sup_{\lambda>0}\frac{h(\lambda x + \bar x,\omega)-h(\bar x,\omega)}{\lambda},
\]
which is independent of the choice of $\bar x\in\dom h(\cdot,\omega)$; see \cite[Theorem~8.5]{roc70a} or \cite[3.21]{rw98}. By \cite[Exercise~14.54(a)]{rw98}, the function $h^\infty$ is a convex normal integrand. If $h(\cdot,\omega)$ has an integrable lower bound, then $h^\infty(x,\omega)\ge 0$ for every $x\in\reals^n$ as is easily seen by letting $\lambda\to\infty$.

\begin{lemma}\label{lem:ip}
Assume that $h_t$ is a normal integrand and that the set-valued mapping
\[
N_t(\omega)=\{x_t\in\reals^{n_t}\,|\,h_t^\infty(x^t,\omega)\le 0,\ x^{t-1}=0\}
\]
is linear-valued. Then $\tilde h_{t-1}$ is a normal integrand with
\[
\tilde h^\infty_{t-1}(x^{t-1},\omega) = \inf_{x_t\in\reals^{n_t}}h^\infty_t(x^{t-1},x_t,\omega).
\]
Moreover, given an $x\in\N$, there is an $\F_t$-measurable $\bar x_t$ such that $\bar x_t(\omega)\perp N_t(\omega)$ and
\[
\tilde h_{t-1}(x^{t-1}(\omega),\omega) = h_t(x^{t-1}(\omega),\bar x_t
(\omega),\omega).
\]
\end{lemma}

\begin{proof}
By \cite[Theorem~9.2]{roc70a}, the linearity condition implies that the infimum in the definition of $\tilde h_{t-1}$ is attained and that $\tilde h_{t-1}(\cdot,\omega)$ is a lower semicontinuous convex function with
\[
\tilde h^\infty_{t-1}(x^{t-1},\omega)=\inf_{x_t\in\reals^{n_t}}h^\infty_t(x^{t-1},x_t,\omega).
\]
By \cite[Proposition~14.47]{rw98}, the lower semicontinuity implies that $\tilde h_{t-1}$ is an $\F_t$-measurable convex normal integrand. By \cite[Proposition~14.45(c)]{rw98}, the function $p(x,\omega):=h_t(x^{t-1}(\omega),x,\omega)$ is then also an $\F_t$-measurable normal integrand so, by \cite[Theorem~14.37]{rw98}, there is an $\F_t$-measurable $\bar x_t$ that attains the minimum for every $\omega$. By \cite[Corollary~8.6.1]{roc70a}, the value of $h_t(x^{t-1}(\omega),x,\omega)$ does not change if we replace $\bar x_t(\omega)$ by its projection to the orthogonal complement of $N_t(\omega)$. By \cite[Exercise~14.17]{rw98}, such a projection preserves measurability.
\end{proof}



It is clear that if $h_t$ has an integrable lower bound, then so will $\tilde h_{t-1}$. Applying Lemmas~\ref{lem:ce} and \ref{lem:ip} recursively backwards for $t=T,\ldots,0$, we then see that if $h$ has an integrable lower bound, the functions $\tilde h_t$ and $h_t$ are well-defined for every $t$ provided that $N_t$ is linear-valued at each step. 

We now get the following refinement of the optimality conditions in \cite[Theorem~1]{rw76} and \cite[Theorems~1 and 2]{evs76} in the convex case.

\begin{theorem}\label{thm:dp}
Assume that $h$ has an integrable lower bound and that $N_t$ is linear-valued for $t=T,\ldots,0$. The functions $h_t$ are then well-defined normal integrands and we have for every $x\in\N$ that
\begin{equation}\label{ie}
Eh_t(x_t(\omega),\omega)\ge \inf\eqref{p}\quad t=0,\ldots,T.
\end{equation}
Optimal solutions $x\in\N$ exist and they are characterized by the condition
\[
x_t(\omega)\in\argmin_{x_t}h_t(x^{t-1}(\omega),x_t,\omega)\quad P\text{-a.s.}\quad t=0,\ldots,T.
\]
which is equivalent to having equalities in \eqref{ie}. Moreover, there is an optimal solution $x\in\N$ such that $x_t\perp N_t$ for every $t=0,\ldots,T$.
\end{theorem}

\begin{proof}
As noted above, a recursive application of Lemmas~\ref{lem:ce} and~\ref{lem:ip} imply that the functions $h_t$ and $\tilde h_t$ are well-defined normal integrands. Given an $x\in\N$, the law of iterated expectations (see e.g.\ Shiryaev~\cite[Section~II.7]{shi96}) gives
\[
Eh_t(x^t(\omega),\omega) \ge E\tilde h_{t-1}(x^{t-1}(\omega),\omega) = Eh_{t-1}(x^{t-1}(\omega),\omega)\quad t=1,\ldots,T.
\]
Thus,
\[
Eh(x(\omega),\omega) = Eh_T(x^T(\omega),\omega) \ge Eh_0(x^0(\omega),\omega) \ge E\inf_{x_0\in\reals^{n_0}}h_0(x_0,\omega),
\]
where the inequalities hold as equalities if and only if
\[
h_t(x^t(\omega),\omega)=\tilde h_{t-1}(x^{t-1}(\omega),\omega)\quad P\text{-a.s.}\quad t=0,\ldots,T.
\]
The existence of such an $x\in\N$ with $x_t\perp N_t$ follows by applying Lemma~\ref{lem:ip} recursively for $t=0,\ldots,T$.
\end{proof}

When the normal integrand $h$ has a separable structure, the dynamic programming equations \eqref{dp} can be written in a more familiar form.

\begin{corollary}[Bellman equations]\label{cor:bel}
Assume that
\[
h(x,\omega) = \sum_{t=0}^T k_t(x_{t-1},x_t,\omega)
\]
for some fixed initial state $x_{-1}$ and $\F_t$-measurable normal integrands $h_t$ with integrable lower bounds. Consider the functions $V_t:\reals^{n_t}\times\Omega\to\ereals$ given by
\begin{equation}\label{bellman}
\begin{split}
V_T(x_T,\omega) &= 0,\\
\tilde V_{t-1}(x_{t-1},\omega) &= \inf_{x_t\in\reals^{n_t}}\{k_t(x_{t-1},x_t,\omega) + V_t(x_t,\omega)\},\\
V_{t-1} &= E_{t-1}\tilde V_{t-1}
\end{split}
\end{equation}
and assume that the set-valued mappings
\[
N_t(\omega) = \{x_t\in\reals^{n_t}\,|\, k_t^\infty(0,x_t,\omega) + V_t^\infty(x_t,\omega)\le 0\}
\]
are linear-valued for each $t=T,\ldots,0$. The functions $V_t$ are then well-defined normal integrands and we have for every $x\in\N$ that
\begin{equation}\label{ieb}
E\left[\sum_{s=0}^tk_s(x_{s-1}(\omega),x_s(\omega),\omega) + V_t(x_t(\omega),\omega)\right] \ge \inf\eqref{p}\quad t=0,\ldots,T.
\end{equation}
Optimal solutions $x\in\N$ exist and they are characterized by the condition
\[
x_t(\omega)\in\argmin_{x_t\in\reals^{n_t}}\{k_t(x_{t-1}(\omega),x_t,\omega) + V_t(x_t,\omega)\}\quad P\text{-a.s.}\quad t=0,\ldots,T,
\]
which is equivalent to having equalities in \eqref{ieb}. Moreover, there is an optimal solution $x\in\N$ such that $x_t\perp N_t$ for every $t=0,\ldots,T$.
\end{corollary}

\begin{proof}
By Theorem~\ref{thm:dp}, it suffices to show that
\begin{equation}\label{Qq}
h_t(x^t,\omega) = \sum_{s=0}^tk_s(x_{s-1},x_s,\omega) + V_t(x_t,\omega)
\end{equation}
for every $t=0,\ldots,T$. For $t=T$, \eqref{Qq} is obvious since $V_T=0$ by definition. Assuming that \eqref{Qq} holds for $t$, we get
\begin{align*}
\tilde h_{t-1}(x^{t-1},\omega) &= \inf_{x_t\in\reals^{n_t}} h_t(x^{t-1},x_t,\omega)\\
&= \sum_{s=0}^{t-1}k_s(x_{s-1},x_s,\omega) + \inf_{x_t\in\reals^{n_t}}\{k_t(x_{t-1},x_t,\omega) + V_t(x_t,\omega)\}\\
&= \sum_{s=0}^{t-1}k_s(x_{s-1},x_s,\omega) + \tilde V_{t-1}(x_{t-1},\omega)
\end{align*}
and then, since for $s=0,\ldots,t-1$, $k_s$ is $\F_{t-1}$-measurable,
\begin{align*}
h_{t-1}(x^{t-1},\omega) &= \sum_{s=0}^{t-1}k_s(x_{s-1},x_s,\omega) + V_{t-1}(x_{t-1},\omega),
\end{align*}
where, by Lemma~\ref{lem:ce}, $V_t$ is a well-defined normal integrand when $N_t$ is linear-valued.
\end{proof}

The rest of this section is devoted to the study of the linearity condition in Theorem~\ref{thm:dp}. Recall that a {\em $\G$-measurable selector} of an $\reals^n$-valued set-valued mapping $C$ is a $\G$-measurable function $x$ such that $x(\omega)\in C(\omega)$ almost surely.

\begin{lemma}\label{lem:rce}
Let $\G\subseteq\F$ be a sigma-algebra and assume that $h$ is an $\F$-normal integrand with an integrable lower bound. If there is an $\bar x\in L^0(\Omega,\G,P;\reals^n)$ such that $Eh(\bar x(\omega),\omega)$ is finite, then $(E^\G h)^\infty = E^\G h^\infty$ and the level sets 
\begin{align*}
\lev_0 h^\infty(\omega) &= \{x\in\reals^n\,|\, h^\infty(x,\omega)\le 0\},\\
\lev_0 (E^\G h)^\infty(\omega) &= \{x\in\reals^n\,|\, (E^\G h)^\infty(x,\omega)\le 0\}
\end{align*}
have the same $\G$-measurable selectors.
\end{lemma}

\begin{proof}
By \cite[Exercise~14.54]{rw98}, $h^\infty$ is a well-defined $\F$-normal integrand. Moreover, the lower bound on $h$ implies that $h^\infty$ is nonnegative. By Lemma~\ref{lem:ce}, $E^\G h$ and $E^\G h^\infty$ are thus well-defined. To show that the latter is the recession function of the former, let $x\in L^0(\Omega,\G,P;\reals^n)$ and $A\in\G$. Convexity of $h$ implies that the difference quotient
\[
\frac{h(\bar x(\omega) +\lambda x(\omega),\omega)-h(\bar x(\omega),\omega)}{\lambda}
\]
is increasing in $\lambda$ for every $\omega$; see e.g.\ \cite[Theorem~23.1]{roc70a}. The lower bound on $h$ and the integrability of $h(\bar x(\cdot),\cdot)$ thus imply that, for $\lambda\ge 1$, the quotients are minorized by a fixed integrable function. Monotone convergence theorem then gives for every $A\in\G$
\begin{align*}
E [1_A h^\infty(x)] &= E [1_A \lim_{\lambda\nearrow\infty} (h(\bar x+\lambda x)-h(\bar x))/\lambda]\\
&=\lim_{\lambda\nearrow\infty} E[1_A (h(\bar x+\lambda x)-h(\bar x))/\lambda]\\
&=\lim_{\lambda\nearrow\infty} E[1_A ((E^\G h)(\bar x+\lambda x)- (E^\G h)(\bar x))/\lambda]\\
&=E[1_A\lim_{\lambda\nearrow\infty} ((E^\G h)(\bar x+\lambda x)- (E^\G h)(\bar x))/\lambda]\\
&=E[1_A (E^\G h)^\infty(x)],
\end{align*}
which means that $(E^\G h)^\infty$ is the conditional expectation of $h^\infty$.

To prove the last claim, let $x\in L^0(\Omega,\G,P;\reals^n)$. By the first claim and the definition of a conditional integrand,
\[
(E^\G h^\infty)(x(\cdot),\cdot)=E^\G h^\infty(x(\cdot),\cdot).
\]
We have $h^\infty(x(\omega),\omega)\le 0$ almost surely if and only if $E^\G h^\infty(x(\cdot),\cdot)\le~0$ almost surely, since $h^\infty\geq 0$.
\end{proof}

\begin{remark}\label{rem:rec}
Consider the parametric problem \eqref{p} and assume that $u\in L^0(\Omega,\F,P;\reals^m)$ is such that $h(\cdot,\omega)=f(\cdot,u(\omega),\omega)$ is proper. We then have
\[
h^\infty(x,\omega)=f^\infty(x,0,\omega),
\]
where $f^\infty(\cdot,\cdot,\omega)$ is the recession function of $f(\cdot,\cdot,\omega)$. It follows that, as soon as they are well defined, the recession functions $\tilde h_t^\infty$ and $h_t^\infty$ and thus, the mappings $N_t$ are independent of the choice of $u\in\dom\varphi$. Indeed, $u\in\dom\varphi$ implies that there is an $x\in\N$ such that $h(x(\omega),\omega)=f(x(\omega),u(\omega),\omega)<\infty$ almost surely. Recursive application of Lemmas~\ref{lem:ip} and \ref{lem:rce} then shows that $h_t$ and $\tilde h_t$ can be expressed in terms of $h^\infty$, which is independent of $u\in\dom\varphi$.
\end{remark}

The following result shows that the linearity condition of Theorem~\ref{thm:dp} can be stated in terms of the original normal integrand $h$ directly. In the proof, we will denote the set of $\G$-measurable selectors of a set-valued mapping $C$ by $L^0(\G;C)$. We will also use the fact that if $C$ is closed-valued and $\G$-measurable, then it is almost surely linear-valued if and only if the set of its measurable selectors is a linear space. This follows easily by considering the Castaing representation of $C$; see e.g.~\cite[Theorem~14.5]{rw98}.

\begin{lemma}\label{lem:lin}
Assume that $h$ has an integrable lower bound and that $E h(\bar x(\omega),\omega)<\infty$ for some $\bar x \in\N$. Then $h_t$ is well-defined and $N_t$ is linear-valued for $t=T,\ldots,0$ if and only if 
\[
\L=\{x\in \N|\, h^\infty(x(\omega),\omega)\le 0 \text{ a.s.}\}
\]
is a linear space. If $x\in\L$ is such that $x^{t-1}=0$ then $x_t\in N_t$ almost surely.
\end{lemma}

\begin{proof}
Redefining $h(x,\omega):=h(x-\bar x(\omega),\omega)$, we may assume that $\bar x=0$. Indeed, such a translation amounts to translating the functions $\tilde h_t$ and $h_t$ accordingly and it does not affect the recession functions $\tilde h_t^\infty$ and $h_t^\infty$. We proceed by induction on $T$. When $T=0$, Lemma~\ref{lem:rce} gives
\begin{align*}
\L&= \{x\in\N|\, h_T^\infty(x(\omega),\omega)\le 0 \text{ a.s.}\} = L^0(\F_T;N_T).
\end{align*}
Since $N_T$ is $\F_T$-measurable, the linearity of $\L$ is equivalent to $N_T$ being linear-valued. Let now $T$ be arbitrary and assume that the claim holds for every $(T-1)$-period model.

If $\L$ is linear then $\L'=\{x\in \N|\, x_0=0,\ h^\infty(x(\omega),\omega)\le 0 \text{ a.s.}\}$ is linear as well. Applying the induction hypothesis to the $(T-1)$-period model obtained by fixing $x_0\equiv0$, we get that $N_t$ is linear for $t=T,\ldots,1$. Applying Lemmas~\ref{lem:ce} and \ref{lem:ip} backwards for $s=T,\ldots,1$, we then see that $h_0$ is well defined. Lemmas~\ref{lem:rce} and \ref{lem:ip} give
\begin{align}
L^0(\F_0;N_0) &= \{x_0\in L^0(\F_0)\,|\, h_0^\infty(x_0(\omega),\omega)\le 0\text{ a.s.}\}\nonumber\\
&= \{x_0\in L^0(\F_0)\,|\, \tilde h_0^\infty(x_0(\omega),\omega)\le 0\text{ a.s.}\}\nonumber\\
&= \{x_0\in L^0(\F_0)\,|\, \inf_{x_1}h_1^\infty(x_0(\omega),x_1,\omega)\le 0\text{ a.s.}\}\nonumber\\
&= \{x_0\in L^0(\F_0)\,|\, \exists\tilde x\in\N:\ \tilde x_0=x_0,\ h_1^\infty(\tilde x^1(\omega),\omega)\le 0\text{ a.s.}\}\nonumber,
\end{align}
where the last equality follows by applying the last part of Lemma~\ref{lem:ip} to the normal integrand $h^\infty$.
Repeating the argument for $t=1,\ldots,T$, we get
\begin{align}
L^0(\F_0;N_0) &= \{x_0\in L^0(\F_0)\,|\, \exists\tilde x\in\N: \tilde x_0=x_0,\ h_T^\infty(\tilde x(\omega),\omega)\le 0 \text{ a.s.}\}\nonumber\\
&= \{x_0\in L^0(\F_0)\,|\, \exists\tilde x\in\N: \tilde x_0=x_0,\ h^\infty(\tilde x(\omega),\omega)\le 0 \text{ a.s.}\}\nonumber\\
&= \{x_0\in L^0(\F_0)\,|\, \exists\tilde x\in\L: \tilde x_0=x_0\}.\label{n0}
\end{align}
The linearity of $\L$ thus implies that of $L^0(\F_0;N_0)$ which is equivalent to $N_0$ being linear-valued.

Assume now that $N_t$ is linear-valued for $t=T,\ldots,0$ and let $x\in\L$. Expression \eqref{n0} for $L^0(\F_0;N_0)$ is again valid so, by linearity of $N_0$, there is an $\tilde x\in\L$ with $\tilde x_0=-x_0$. Since $h^\infty$ is sublinear, $\L$ is a cone, so that $x+\tilde x\in\L$. Since $x_0+\tilde x_0=0$, we also have $x+\tilde x\in\L'$. Since, by the induction assumption, $\L'$ is linear and since $\L'\subseteq\L$, we get $-x-\tilde x\in\L$. Since $\L$ is a cone, we get $-x=\tilde x-x-\tilde x\in\L$. Thus, $\L$ is linear.

For $t=0$, the last claim follows directly from expression \eqref{n0}. The general case follows by applying this to the $(T-t)$-period model obtained by fixing $x^{t-1}\equiv 0$.
\end{proof}

When $h$ is the indicator function of a convex set, the linearity condition in Lemma~\ref{lem:lin} becomes the linearity condition of \cite[Theorem~8]{pen11c} which generalizes various no-arbitrage conditions that have been used in mathematical finance. The following example illustrates the situation in the classical perfectly liquid market model; see \cite{pen11c} for more general models.

\begin{example}[Superhedging in liquid markets]\label{ex:sh}
Let $S=(S_t)_{t=0}^T$ be an $\reals^d$-valued $\FF$-adapted stochastic process, $n_t=d$, $m=1$ and
\[
f(x,u,\omega) = 
\begin{cases}
0 & \text{if $\sum_{t=0}^{T-1}x_t\cdot\Delta S_t(\omega)\ge u$,}\\
+\infty & \text{otherwise}.
\end{cases}
\]
We get
\[
\varphi(u)=\inf_{x\in\N} I_f(x,u) =
\begin{cases}
0 & \text{if $×u\in\C$},\\
+\infty & \text{otherwise},
\end{cases}
\]
where $\C = \{u\in L^0\,|\,\exists x\in\N:\ \sum_{t=0}^{T-1}x_t\cdot\Delta S_t\ge u\}$. In the classical perfectly liquid model of financial markets, where $S$ gives the unit prices of the ``risky assets'' and $x_t$ is the portfolio held over $(t,t+1]$, the set $\C$ consists of the {\em contingent claims} that can be {\em superhedged} without a cost; see e.g.~\cite[Section~6.4]{ds6}. Since $f$ is a closed positively homogeneous function, we have $f^\infty=f$ and
\[
\{x\in \N\,|\, f^\infty(x(\omega),0,\omega)\le 0 \text{ a.s.}\}=\{x\in \N|\, \sum_{t=0}^{T-1}x_t\cdot\Delta S_t\ge 0\}.
\]
This set is linear, and thus, the function $h(x,\omega)=f(x,u(\omega),\omega)$ satisfies the linearity condition in Lemma~\ref{lem:lin}, if and only if the price process $S$ satisfies the {\em no-arbitrage} condition. 
\end{example}

The following simple example goes beyond indicator functions and also of inf-compact integrands considered in \cite{rw76,evs76}.

\begin{example}[Variance optimal hedging]\label{ex:voh}
Let $S=(S_t)_{t=0}^T$ be an $\reals^d$-valued $\FF$-adapted stochastic process, $u\in L^0(\Omega,\F,P;\reals)$ and consider the problem of minimizing
\[
E(V_0 + \sum_{t=0}^{T-1} z_t\cdot\Delta S_{t+1} - u)^2
\]
over $V_0\in\reals$ and $\F_t$-measurable $\reals^d$-valued functions $z_t$. This corresponds to \eqref{p} with $x_0=(z_0,V_0)$, $x_t=z_t$ for $t=1,\ldots,T$ and 
\[
f(x,u,\omega) = (V_0 + \sum_{t=0}^{T-1} z_t\cdot\Delta S_{t+1}(\omega) - u)^2.
\]
The above problem has been studied e.g.\ in F\"ollmer and Schied~\cite[Section~10.3]{fs4}, where $V_0$ is interpreted as an initial value of a self-financing trading strategy where $z_t$ is the portfolio of risky assets held over period $[t,t+1]$. By \cite[Theorem~9.4]{roc70a},
\[
f^\infty(x,u,\omega) = 
\begin{cases}
0 & \text{if $V_0 + \sum_{t=0}^{T-1} z_t\cdot\Delta S_{t+1}(\omega) - u=0$},\\
+\infty & \text{otherwise}.
\end{cases}
\]
By Remark~\ref{rem:rec}, the function $h(x,\omega)=f(x,u(\omega),\omega)$ then satisfies the linearity condition of Lemma~\ref{lem:lin}, so the optimal solution is attained. This should be compared with the existence results in \cite[Section~10.3]{fs4}, where it was assumed that $d=1$.
\end{example}

\section{Lower semicontinuity of the value function}\label{sec:lsc}

We now return to the parametrized problem \eqref{p}. Being the inf-projection of the convex integral functional
\[
I_f(x,u)=Ef(x(\omega),u(\omega),\omega),
\]
the value function 
\[
\varphi(u)=\inf_{x\in\N}Ef(x(\omega),u(\omega),\omega)
\]
is convex on $L^0(\Omega,F,P;\reals^m)$; see e.g.\ \cite[Theorem~1]{roc74}. Our aim is to give conditions under which $\varphi$ is lower semicontinuous on certain locally convex topological vector subspaces of $L^0(\Omega,F,P;\reals^m)$. The lower semicontinuity is equivalent to the absence of a duality gap in the duality framework of \cite{pen11c} (which is essentially an instance of the conjugate duality framework of Rockafellar~\cite{roc74}) which we now briefly recall (and slightly generalize).

Assume that $\U$ and $\Y$ are vector subspaces of $L^0(\Omega,F,P;\reals^m)$ in {\em separating duality} under the bilinear form
\[
\langle u,y\rangle = E[u(\omega)\cdot y(\omega)],
\]
i.e.\ that $E[u(\omega)\cdot y(\omega)]$ is finite for every $u\in\U$ and $y\in\Y$ and that for every nonzero $u\in\U$ (resp.\ $y\in\Y$), there is at least one $y\in\Y$ (resp.\ $u\in\U$) such that $\langle u,y\rangle\ne 0$. The special case $\U=L^p$ and $\Y=L^q$ was studied in \cite{pen11c}. The weakest and the strongest locally convex topologies on $\U$ compatible with the pairing will be denoted by $\sigma(\U,\Y)$ and $\tau(\U,\Y)$, respectively. Since the value function $\varphi$ is convex, we have by the classical separation argument, that $\varphi$ is lower semicontinuous with respect to $\sigma(\U,\Y)$ if it is merely lower semicontinuous with respect to $\tau(\U,\Y)$. When $\U=L^p$ and $\Y=L^q$ for $p\in[1,\infty)$, $\tau(\U,\Y)$ is simply the norm topology on $\U$ and $\sigma(\U,\Y)$ the weak topology. A general treatment of topological spaces in separating duality can be found e.g.\ in Kelley and Namioka~\cite{kn76}.

The {\em Lagrangian} associated with \eqref{p} is the extended real-valued function
\[
L(x,y) = \inf_{u\in\U}\{I_f(x,u)-\langle u,y\rangle\}
\]
on $\N\times\Y$. The Lagrangian is convex in $x$ and concave in $y$. The {\em dual objective} is the extended real-valued function on $\Y$ defined by 
\[
g(y) = \inf_{x\in\N}L(x,y).
\]
The basic duality result \cite[Theorem~7]{roc74} says, in particular, that $g = -\varphi^*$. When $\varphi$ is lower semicontinuous and proper, the biconjugate theorem (see e.g.\ \cite[Theorem~5]{roc74}) then gives the dual representation
\begin{equation}\label{bc}
\varphi(u) = \sup\{\langle u,y\rangle + g(y)\}.
\end{equation}
It was shown in \cite{pen11c} that this abstract result is behind many duality frameworks in stochastic optimization and mathematical finance.

It was assumed in \cite{pen11c} that $\U=L^p$ and $\Y=L^q$, but the main result \cite[Theorem~3]{pen11c} remains valid as long as the space $\U$ is {\em decomposable} in the sense that
\[
\chi_Au+\chi_{\Omega\setminus A}u'\in\U
\]
whenever $A\in\F$, $u\in\U$ and $u'\in L^\infty(\Omega,\F,P;\reals^m)$. Indeed, the decomposability property allows the use of the interchange rule for minimization and integration (see \cite[Theorem~14.60]{rw98}), which suffices for the proof of \cite[Theorem~3]{pen11c}. We also note that decomposability of the spaces $\U$ and $\Y$ implies that the separation property holds automatically for the bilinear form defined above; see \cite[Lemma~6]{val75}. 
Moreover, we have the following relations for relative topologies.

\begin{lemma}\label{lem:rel}
If $\U$ and $\Y$ are decomposable, then $L^\infty\subseteq\U\subseteq L^1$ and
\begin{align*}
\sigma(L^1,L^\infty)|_\U&\subseteq\sigma(\U,\Y),\quad\sigma(\U,\Y)|_{L^\infty}\subseteq\sigma(L^\infty,L^1),\\
\tau(L^1,L^\infty)|_\U&\subseteq\tau(\U,\Y),\quad\tau(\U,\Y)|_{L^\infty}\subseteq\tau(L^\infty,L^1).
\end{align*}
\end{lemma}

\begin{proof}
By \cite[Lemme 1, p.??]{bis73}, $L^\infty\subseteq\U\subset L^1$ and $L^\infty\subseteq\Y\subseteq L^1$ which give the relations for the $\sigma$-topologies. Since, by symmetry, analogous relations are valid for the $\sigma$-topologies on $\Y$, we have that $\sigma(L^\infty,L^1)$-compact subsets of $L^\infty$ are $\sigma(\Y,\U)$-compact. Since, by the Mackey-Arens theorem, $\tau(\U,\Y)$ is generated by the support functions of $\sigma(\Y,\U)$-compact sets, we get $\tau(L^1,L^\infty)|_\U\subseteq\tau(\U,\Y)$. The remaining inclusion is verified similarly.
\end{proof}

The traditional ``direct method'' for proving the lower semicontinuity would be to assume that the integral functional $I_f$ is uniformly inf-compact in $x$ with respect to an appropriate topology on $\N$. If the topology is strong enough to imply the almost sure convergence of a subsequence, the sequential lower semicontinuity can often be derived from Fatou's lemma. In certain applications, this purely topological argument fails because $I_f$ lacks an appropriate inf-compactness property in $x$. In convex problems, the following measure theoretic result can sometimes be used as a substitute for compactness.

\begin{lemma}[Koml\'os' theorem]\label{lem:kom}
Let $(x^\nu)_{\nu=1}^\infty$ be a sequence in $L^0(\Omega,\F,P;\reals^n)$ which is either
\begin{enumerate}
\item
bounded in $L^1$,
\item
almost surely bounded in the sense that 
\[
\sup_\nu|x^\nu(\omega)|<\infty\quad P\text{-a.s.}
\]
\end{enumerate}
Then there is a sequence of convex combinations $\bar x^\nu\in\co\{x^\mu\,|\,\mu\ge\nu\}$ that converges almost surely to an $\reals^n$-valued function.
\end{lemma}

\begin{proof}
See e.g.\ \cite{ds6} or \cite{ks9}.
\end{proof}

The following is our main result.

\begin{theorem}\label{thm:lsc}
Assume that there is a $y\in\Y$ and an $m\in L^1(\Omega,\F,P)$ such that for $P$-almost every $\omega$,
\[
f(x,u,\omega)\ge u\cdot y(\omega) +m(\omega)\quad \forall (x,u)\in\reals^n\times\reals^m 
\]
and that $\{x\in \N|\, f^\infty(x(\omega),0,\omega)\le 0 \text{ a.s.}\}$ is a linear space. Then
\[
\varphi(u) = \inf_{x\in\N}I_f(x,u)
\]
is lower semicontinuous on $\U$ and the infimum is attained for every $u\in\U$.
\end{theorem}

\begin{proof}
Let $h_u(x,\omega)=f(x,u(\omega),\omega)$. The lower bound on $f$ implies that $h_u$ has an integrable lower bound. As noted in Remark~\ref{rem:rec}, $h_u^\infty(x,\omega)=f^\infty(x,0,\omega)$ for every $u$, so the linearity condition on $f$ implies that $h_u$ satisfies the linearity conditions in Lemma~\ref{lem:lin}. By Theorem~\ref{thm:dp}, the infimum in $\varphi(u)= \inf_{x\in\N}I_f(x,u)$ is thus attained for every $u\in\U$ by an $x\in\N$ with $x_t(\omega)\perp N_t(\omega)$ almost surely.

For lower semicontinuity, it suffices to show that $\varphi$ is lower semicontinuous on the linear space $L^{1,y}=\{u\in L^1|\, |E[u\cdot y]|<\infty\}$ with respect to the norm $|| u||_{L^{1,y}}=E|u|+ |E[u\cdot y]|$. Indeed, since $y\in\Y$, we have $\U\subseteq L^{1,y}$ and, by Lemma~\ref{lem:rel}, the norm $||\cdot||_{L^{1,y}}$ is continuous on $\tau(\U,\Y)$, which means that $\tau(\U,\Y)$ is stronger than the norm topology restricted to $\U$. Since $L^{1,y}$ is a normed space, it suffices to prove sequential lower semicontinuity, which means that for any $\gamma\in\reals$ and for any sequence $(u^\nu)_{\nu=1}^\infty$ such that 
\[
\varphi(u^\nu)\le\gamma
\]
and $u^\nu\to u$ in $L^{1,y}$, we have $\varphi(u)\le\gamma$. We will prove this by establishing the existence of an $x\in\N$ such that $I_f(x,u)\le\gamma$.

As observed at the beginning of the proof, there is for every $\nu$ an $x^\nu\in\N$ such that $x^\nu_t\perp N_t$ and
\[
I_f(x^\nu,u^\nu)\le \gamma.
\]
Moreover, the mappings $N_t$ are independent of $u^\nu$; see Remark~\ref{rem:rec}.
Since $u^\nu$ converges in $L^{1,y}$, the lower bound on $f$ implies that the negative parts of the functions $\omega\mapsto f(x^\nu(\omega),u^\nu(\omega),\omega)$ are bounded in $L^1$. Since $I_f(x^\nu,u^\nu)\le \gamma$, the positive parts must be bounded as well. Thus, by Lemma~\ref{lem:kom}, there is a sequence of convex combinations 
\[
\phi^\nu(\omega):=\sum_{\mu=\nu}^\infty\alpha^{\nu,\mu}f(x^\mu(\omega),u^\mu(\omega),\omega)
\]
that converges almost surely to a real-valued measurable function. In particular, the function $\phi(\omega):=\sup_\nu\phi^\nu(\omega)$ is almost surely finite. Defining 
\[
(\bar x^\nu,\bar u^\nu) = \sum_{\mu=\nu}^\infty \alpha^{\nu,\mu}(x^\mu,u^\mu)
\]
we have by convexity that 
\[
f(\bar x^\nu(\omega),\bar u^\nu(\omega),\omega)\le\phi^\nu(\omega)\le\phi(\omega)\quad P\text{-a.s.}
\]
and $I_f(\bar x^\nu,\bar u^\nu)\le \gamma$. Moreover, we still have $\bar x^\nu_t\in N_t^\perp$ almost surely and $\bar u^\nu\to u$ in the $L^{1,y}$-norm. 

Passing to a subsequence if necessary, we may assume that $\bar u^\nu\to u$ almost surely, so that the measurable function $\rho(\omega) := \sup_\nu|\bar u^\nu(\omega)|$ is almost surely finite. Each $(\bar x^\nu,\bar u^\nu)$ then belongs to the set
\[
\C=\{(x,u)\in\N\times L^0\,|\, (x,u)\in C\ \text{a.s.}\},
\]
where $C(\omega) = \{(x,u)\,|\, x_t\in N_t^\perp(\omega),\ u\in\rho(\omega)\uball,\ f(x,u,\omega)\le\phi(\omega)\}$. We will now apply \cite[Theorem~6]{pen11c}, which says that the sequence $(\bar x^\nu,\bar u^\nu)_{\nu=1}^\infty$ is almost surely bounded if  
\begin{equation}\label{zero}
\{(x,u)\in\N\times L^0\,|\, (x,u)\in C^\infty\ \text{a.s.}\}=\{(0,0)\}.
\end{equation}
By Corollary~8.3.3 and Theorem~8.7 of \cite{roc70a}, 
\[
C^\infty(\omega)=\{(x,0)\,|\, x_t\in N_t^\perp(\omega),\ f^\infty(x,0,\omega)\le 0\}.
\]
If $x\in\N$ is such that $f^\infty(x(\omega),0,\omega)\le 0$ then, by the last part of Lemma~\ref{lem:lin}, we have $x_0\in N_0$. The condition $x_0\in N_0^\perp$ then implies that $x_0=0$. Repeating the argument for $t=1,\ldots,T$ gives \eqref{zero} so $(\bar x^\nu,\bar u^\nu)_{\nu=1}^\infty$ is almost surely bounded.

By Lemma~\ref{lem:kom}, there is a sequence $(\hat x^\nu,\hat u^\nu)_{\nu=1}^\infty$ of convex combinations of $(\bar x^\nu,\bar u^\nu)_{\nu=1}^\infty$ that converges almost surely to a point $(x,\hat u)$, where necessarily $\hat u=u$ since $\bar u^\nu\to u$ almost surely. We still have $\hat u^\nu\to u$ in the $L^{1,y}$-norm and, by convexity, $I_f(\hat x^\nu,\hat u^\nu)\le\gamma$. By Fatou's lemma,
\begin{multline*}
E[f(x(\omega),u(\omega),\omega)-y(\omega)\cdot u(\omega)-m(\omega)]\\
\le\liminf_{\nu\to\infty} E[f(\hat x^\nu(\omega),\hat u^\nu(\omega),\omega)-y(\omega)\cdot\hat u^\nu(\omega)-m(\omega)],
\end{multline*}
where $E[y(\omega)\cdot\hat u^\nu(\omega)]\to E[y(\omega)\cdot u(\omega)]$, by the $L^{1,y}$-convergence, so that
\[
I_f(x,
u)\le \liminf_{\nu\to\infty}I_f(\hat x^\nu,\hat u^\nu)\le \gamma,
\]
which completes the proof.
\end{proof}

\section{An application to mathematical finance}\label{sec:app}

We will illustrate Theorem~\ref{thm:lsc} on the optimal consumption problem considered in \cite[Section~5]{pen11c}. The problem is set in a generalization of the market model of Kabanov~\cite{kab99}, where a finite number $d$ of securities is traded over finite discrete time $t=0,\ldots,T$. At each time $t$ and state $\omega\in\Omega$, the market is described by two closed convex sets, $C_t(\omega)$ and $D_t(\omega)$, both of which contain the origin. The set $C_t(\omega)$ consists of the portfolios that are freely available in the market and $D_t(\omega)$ consists of the portfolios that the investor is allowed to hold over the period $[t, t+1)$. For each $t$, the sets $C_t$ and $D_t$ are assumed to be $\F_t$-measurable. 

Consider the problem
\begin{align}\label{ocp}
\begin{split}
\maximize_{(z,c)\in\N}\quad & E\sum_{t=0}^T U_t(c_t)\\
\st \quad  & z_t-z_{t-1} + c_t \in C_t,\ z_t\in D_t\quad P\text{-a.s.}\ t=0,\ldots,T,
\end{split}
\end{align}
where $z_{-1}:=0$, $D_T(\omega):=\{0\}$ and $-U_t$ is a convex $\F_t$-measurable normal integrand on $\reals^d\times\Omega$. This models an optimal consumption problem where at each time $t$ and stage $\omega$ we can consume some of the assets and update the existing portfolio $z_{t-1}$. The combined process $(z,c)$ is required to be {\em self-financing} in the sense that the sum of the portfolio update $\Delta z_t:=z_t-z_{t-1}$ and the consumption vector $c_t$ has to be freely available in the market, i.e.\ it belongs to $C_t(\omega)$. In addition, the portfolio constraint $z_t(\omega)\in D_t(\omega)$ is required to hold almost surely at each time. Problem \eqref{ocp} generalizes the classical optimal consumption problem where the numeraire asset is consumed in a perfectly liquid market model (see Examples~\ref{ex:sh} and \ref{ex:voh}). A general treatment of the continuous-time model can be found in Karatzas and \u{Z}itkovi\'c~\cite{kz3}.

Defining
\[
\C = \{c\in\A\,|\, \Delta z_t + c_t \in C_t,\ z_t\in D_t\quad P\text{-a.s.}\ t=0,\ldots,T\},
\]
where $\A$ denotes the set of $\reals^d$-valued adapted processes (so that $\N=\A\times\A$), we can write problem \eqref{ocp} compactly as
\begin{equation}\label{ocpc}
\maximize_{c\in\A}\quad E\sum_{t=0}^T U_t(c_t)\quad\text{over $c\in\C$}.
\end{equation}
The set $\C$ can be interpreted as the set of consumption processes that can be super-replicated without a cost in the market given by the pair $(C,D)$; compare with the definition of the set $\C$ in Example~\ref{ex:sh}.

In order to dualize the problem, we embed it in the general duality framework with, $x_t=(z_t,c_t)$, $u=(u_t)_{t=0}^T$ and 
\begin{align*}
f(x,u,\omega)= 
\begin{cases}
-\sum_{t=0}^T U_t(c_t,\omega) & \text{if }\Delta z_t + c_t +u_t \in C_t(\omega),\, z_t\in D_t(\omega)\\
+\infty &\text{otherwise}.
\end{cases}
\end{align*}
Here $u_t\in\reals^d$ so that the dimension of $u$ equals $m=(T+1)d$. The Lagrangian integrand becomes
\begin{align*}
&\ l(x,y,\omega) = \inf_{u\in\reals^m}\{f(x,u,\omega)-u\cdot y\}\\
&=
\inf_{u\in\reals^m}\{-\sum_{t=0}^T[U_t(c_t,\omega)+u_t\cdot y_t]\, | \, \Delta z_t + c_t +u_t \in C_t(\omega),\ z_t\in D_t(\omega)\}\\
&=
\inf_{\tilde u\in\reals^m}\{-\sum_{t=0}^T [U_t(c_t,\omega)+(\tilde u_t-\Delta z_t-c_t)\cdot y_t]\, | \, \tilde u_t \in C_t(\omega),\ z_t\in D_t(\omega)\}\\
&= 
\begin{cases}
-\sum_{t=0}^T [U_t(c_t,\omega)+\sigma_{C_t(\omega)}(y_t) - (\Delta z_t+c_t)\cdot y_t] & \text{if $z_t\in D_t(\omega)$}\\
+\infty & \text{otherwise}
\end{cases}\\
&= 
\begin{cases}
-\sum_{t=0}^T [U_t(c_t,\omega) + \sigma_{C_t(\omega)}(y_t) + z_t\cdot\Delta y_{t+1} - c_t\cdot y_t] & \text{if $z_t\in D_t(\omega)$}\\
+\infty & \text{otherwise},
\end{cases}
\end{align*}
where $\sigma_{C_t(\omega)}$ denotes the support function of $C_t(\omega)$. In the last equality we have used the ``integration by parts'' formula
\[
\sum_{t=0}^T\Delta z_t\cdot y_t = -\sum_{t=0}^Tz_t\cdot\Delta y_{t+1}
\]
where $y_{T+1}:=0$. 

Recall that, by Lemma~\ref{lem:rel}, $\Y\subset L^1$. On the other hand, by \cite[Theorem 3]{pen11c}\footnote{Theorem 3 of \cite{pen11c} is stated for the case $\U=L^p$ and $\Y=L^q$, but its proof goes through in the general case without a change.},
\begin{align*}
g(y)=\inf_{x\in\N^\infty} El(x(\omega),y(\omega),\omega),
\end{align*}
where $\N^\infty=\N\cap L^\infty$. We can then use the law of iterated expectations (see e.g.\ \cite[Section~II.7]{shi96}) and the interchange rule for integration and minimization (see e.g.\ \cite[Theorem~14.60]{rw98}) to write the dual objective as
\begin{align*}
g(y) &= \inf_{(c,z)\in\N^\infty}\left\{ E\sum_{t=0}^T[-U_t(c_t)-\sigma_{C_t}(y_t) \right.\\
&\left. \qquad\qquad\qquad\qquad {} - z_t\cdot E_{t}\Delta y_{t+1}+c_t\cdot E_t y_t] \left.\right|\,  z_t\in D_t\text{ a.s.} \vphantom{E\sum_{t=0}^T}\right\}\\
&= E\sum_{t=0}^T \inf_{c_t,z_t\in\reals^d}\left\{-U_t(c_t,\omega)-\sigma_{C_t(\omega)}(y_t(\omega)) \right.\\
 &\left. \qquad\qquad\qquad\quad {}-z_t\cdot(E_t\Delta y_{t+1})(\omega)+c_t\cdot(E_ty_t)(\omega)|\, z_t\in D_t(\omega)\right\}\\
&= E\sum_{t=0}^T [U^*_t(E_ty_t) -\sigma_{C_t}(y_t)-\sigma_{D_t}(E_t\Delta y_{t+1})].
\end{align*}
where 
\[
U_t^*(y,\omega) = \inf_{c\in\reals^d}\{c\cdot y-U_t(c,\omega)\}
\]
is the conjugate of $U_t$ in the concave sense.

When $C_t(\omega)$ and $D_t(\omega)$ are convex cones, we have
\[
g(y)=
\begin{cases}
E\sum_{t=0}^T U^*_t(E_ty_t) & \text{if $y\in\D$},\\
-\infty & \text{otherwise},
\end{cases}
\]
where
\[
\D = \{y\in\Y\,|\,y_t\in C_t^*,\ E_t\Delta y_t\in D_t^*\},
\]
where $C_t^*(\omega)$ and $D_t^*(\omega)$ are the polar cones of $C_t(\omega)$ and $D_t(\omega)$, respectively. The dual problem can then be written as
\begin{equation}\label{ocpd}
\maximize\quad E\sum_{t=0}^TU^*_t(y_t)\quad\text{{\rm over}}\ y\in\D
\end{equation}
in symmetry with the primal problem \eqref{ocpc}. The $\FF$-adapted elements of $\D$ are called {\em consistent price systems} for the market model $(C,D)$; see \cite[Example~4.2]{pen11c}. It follows from Jensen's inequality that the value of the dual objective does not decrease when replacing a general $y\in\Y$ by its $\FF$-adapted projection; see \cite[Example~3.6]{pen11c}.

In \cite{pen11c}, the lower semicontinuity of the value function associated with the optimal consumption problem was left open so it could not be claimed that the optimal values of \eqref{ocpc} and \eqref{ocpd} are equal. With the help of Theorem~\ref{thm:lsc}, we can now derive simple sufficient conditions. We will assume that the utility function $U_t$ satisfies the growth condition
\begin{equation}\label{gc}
U_t^\infty(c,\omega)=\begin{cases}
0 & \text{if $c\in\reals_+^d$},\\
-\infty & \text{otherwise},
\end{cases}
\end{equation}
for every $\omega$. When $d=1$ and $U_t(\cdot,\omega)$ is smooth, this is equivalent to the conditions $\lim_{c\to\infty}U_t'(c,\omega)=0$ and $\lim_{c\to-\infty}U_t'(c,\omega)=+\infty$ which generalize the Inada conditions; see e.g.\ \cite{ks99}. We will say that the market model $(C,D)$ satisfies the condition of {\em no-scalable arbitrage} if 
\begin{equation}\label{nsa}
\{c\in\A_+\,|\,\exists z\in\A:\ \Delta z_t + c_t\in C^\infty_t,\ z_t\in D^\infty_t\}=\{0\},
\end{equation}
where $\A_+$ denotes the set of componentwise nonnegative claim processes. This condition is related to arbitrage opportunities that can be scaled up by arbitrarily large positive numbers; see \cite{pen11,pp10}.

\begin{theorem}\label{thm:ocp}
Assume that the optimal value of \eqref{ocp} or, equivalently, \eqref{ocpc} is less than $\infty$, that $U_t$ satisfy the growth condition \eqref{gc} and that there is an integrable function $m$ such that $U_t(c,\omega)\le m(\omega)$ almost surely for every $c\in\reals^d$ and $t=0,\ldots,T$. If the market model $(C,D)$ satisfies the no-scalable arbitrage condition \eqref{nsa} and if the set
\begin{align*}
&\{z\in\A\,|\, \Delta z_t \in C^\infty_t,\ z_t\in D^\infty_t\},
\end{align*}
is linear, then the optimal value of
\[
\maximize\quad E\sum_{t=0}^T U_t(c_t)\quad\text{over}\quad c\in\C-u
\]
is lower semicontinuous as a function of $u\in\U$ and the infimum is always attained. In particular, if $C$ and $D$ are conical, then the optimal value of the primal \eqref{ocpc} is the negative of the optimal value of the dual \eqref{ocpd}.
\end{theorem}

\begin{proof}
By \cite[Theorems~9.3 and 9.5]{roc70a},
\begin{align*}
f^\infty(x,u,\omega)=
\begin{cases}
-\sum_{t=0}^T U^\infty_t(c_t,\omega) & \text{if }\Delta z_t + c_t \in C^\infty_t(\omega),\, z_t\in D^\infty_t(\omega)\\
+\infty & \text{otherwise},
\end{cases}
\end{align*}
so, by Lemma~\ref{lem:lin}, the linearity condition of Theorem~\ref{thm:lsc} means that the set
\[
\L=\{(z,c)\in\N \, | \, \sum_{t=0}^TU_t^\infty(c_t)\ge 0,\ \Delta z_t + c_t\in C^\infty_t,\, z_t\in D^\infty_t\}
\]
is linear. Under the growth condition \eqref{gc}, 
\[
\L=\{(z,c)\in\N \, | \, c_t\ge 0,\ \Delta z_t + c_t\in C^\infty_t,\, z_t\in D^\infty_t\}
\]
If $(z,c)\in\L$, condition \eqref{nsa} implies that $c=0$, and then, by linearity of the set $\{z\in\A\,|\, \Delta z_t \in C^\infty_t,\ z_t\in D^\infty_t\}$, we have $(-z,-c)\in\L$. Since $\L$ is also a cone, it has to be a linear space.
\end{proof}

\begin{remark}\label{rem:rnsa}
The conclusions of Theorem~\ref{thm:ocp} remain valid if, instead of the growth condition \eqref{gc} and the no-scalable arbitrage condition \eqref{nsa}, we assume that the set
\begin{equation}\label{lin3}
\{c\in\A\,|\,\exists z\in\A:\ \sum_{t=0}^TU_t^\infty(c_t)\ge 0,\ \Delta z_t + c_t\in C^\infty_t,\ z_t\in D^\infty_t\},
\end{equation}
is linear. Indeed, if $\L$ is as in the above proof and $(z,c)\in\L$, condition \eqref{lin3} gives the existence of a $z^-\in\A$ such that $(z^-,-c)\in\L$. Since $\L$ is a cone, we get $(z+z^-,0)\in\L$ and then, the linearity condition of Theorem~\ref{thm:ocp} gives $(-(z+z^-),0)\in\L$. Since $-(z,c)=(-(z+z^-),0)+(z^-,-c)$, we get that $-(z,c)\in\L$, i.e., $\L$ is linear.
\end{remark}

\bibliographystyle{plain}
\bibliography{sp}

\end{document}